\documentclass[11pt]{article}

\usepackage[T2A]{fontenc}

\usepackage[cp1251]{inputenc}

\usepackage{amsfonts}

\usepackage{eufrak}

\usepackage{amssymb}

\usepackage{amsmath}


\begin{document}

\noindent{\textbf{Characterization theorems for $Q$-independent random}}

\noindent{\textbf{variables with values in a locally compact Abelian group}}

\bigskip

\noindent {Gennadiy Feldman}

\bigskip


\noindent{\small{{\bf Abstract.} Let  $X$ be a locally compact Abelian group, $Y$ be its character group.  Following A. Kagan and G. Sz\'ekely we introduce
 a notion of  $Q$-independence for random variables with values in $X$.
We prove group analogues of the Cram\'er, Kac--Bernstein, Skitovich--Darmois and Heyde theorems for $Q$-independent random variables with values in $X$. The proofs of these theorems are reduced to solving some functional equations
on the group $Y$.}}

\bigskip


\bigskip


\bigskip

\noindent{\textbf{1. Introduction}}

\bigskip

\noindent It is well known that if a random variable  $\xi$ has a Gaussian distribution
and $\xi$ is a sum of two independent random variables
$\xi=\xi_1+\xi_2$, then $\xi_j$ are also Gaussian
 (Cram\'er's theorem). Let $\xi_1, \dots, \xi_n$ be independent random variables. Consider the linear forms
$L_1=\alpha_1\xi_1+\cdots+\alpha_n\xi_n$ and
$L_2=\beta_1\xi_1+\cdots+\beta_n\xi_n$, where the coefficients
$\alpha_j,\beta_j$ are nonzero real numbers. Then the independence
of the linear forms  $L_1$ and $L_2$
implies that the random variables   $\xi_j$ are Gaussian (Skitovich--Darmois's theorem). A theorem similar
to the Skitovich--Darmois theorem was proved by
Heyde. In this theorem a Gaussian distribution
is characterized by the symmetry of the conditional distribution
of the linear form  $L_2$ given
$L_1$. On the one hand, in the article  \cite{KS} A. Kagan and G. Sz\'ekely introduced
a notion of $Q$-independence and proved, in particular,
that the classical Cram\'er and   Skitovich--Darmois theorems
hold true if instead of independence  $Q$-independence is considered.
On the other hand, in the  articles  \cite{Fe1}, \cite{Fe3}, \cite{Fe16}, \cite{Fe30}, see also
\cite{Fe5a},
the locally compact Abelian groups $X$ for which
group analogues of the Cram\'er,   Skitovich--Darmois and Heyde theorems
for independent random variables with values in  $X$ where described.

In this article we introduce the notion of $Q$-independence for random variables taking values in
a locally compact Abelian group. We prove that if we consider  $Q$-independence instead of independence, then group analogues of theorems by Cram\'er,   Skitovich--Darmois and Heyde
hold true for the same classes of groups. We also consider a group analogue  of the
Kac--Bernstein theorem for $Q$-independent random variables.
The proofs of the corresponding theorems are reduced to solving some functional equations
 on the character group of the initial group in the class of continuous positive definite functions. Standard results on abstract harmonic analysis (see e.g. \cite{HeRo1})
will be used.

 In this paper we suppose $X$ is a second countable locally compact Abelian group. Denote by $Y$ the character group of $X$, and by $(x,y)$ the value of a character $y \in Y$ at $x\in X$.
  If $\xi$ is a  random variable with values in the group
 $X$, then denote by $\mu_\xi$ its distribution and by
$$
\hat\mu_\xi(y)={\bf E}[(\xi, y)]=\int_{X}(x, y)d \mu_\xi(x),
\quad y\in Y,$$   the characteristic function of the
distribution $\mu_\xi$. We will also call  $\hat\mu_\xi(y)$
the characteristic function of the random variable $\xi$.
Let $f(y)$ be a function on the group    $Y$,   and let $h \in
Y$. Denote by   $\Delta_h$   the finite difference operator
$$
\Delta_h f(y)=f(y+h)-f(y).
$$
A function $f(y)$ on $Y$ is called
a    polynomial  if
$$\Delta_{h}^{n+1}f(y)=0$$ for some   $n$ and for all $y, h \in Y$.
The minimal $n$ for which this equality holds is called   the
degree  of the polynomial $f(y)$.

Let $\xi_1, \dots, \xi_n$ be random variables with values in the group
 $X$.
Following A. Kagan and G. Sz\'ekely  (\cite{KS}) we say that the random variables
  $\xi_1,  \dots, \xi_n$ are
$Q$-independent  if their join characteristic function is represented in the form
\begin{equation}\label{17}
\hat\mu_{(\xi_1, \dots, \xi_n)}(y_1, \dots, y_n)={\bf E}[(\xi_1, y_1)\cdots(\xi_n, y_n)]=$$$$=\left(\prod_{j=1}^n\hat\mu_{\xi_j}(y_j)\right)\exp\{q(y_1,
\dots, y_n)\}, \quad y_j\in Y,
\end{equation}
where $q(y_1, \dots, y_n)$ is a continuous polynomial on the group
 $Y^n$.
We will also assume that $q(0, \dots, 0)=0$.

Denote by ${\rm M}^1(X)$   the convolution semigroup of probability distributions on
the group   $X$.  We remind that
a distribution  $\gamma\in {\rm M}^1(X)$  is called Gaussian
(see \cite[Chapter IV]{Pa}),
if its characteristic function is represented in the form
\begin{equation}\label{1}
\hat\gamma(y)=(x,y)\exp\{-\varphi(y)\}, \quad y\in Y,
\end{equation}
where $x \in X$, and $\varphi(y)$ is a continuous non-negative function
on the group $Y$
 satisfying the equation
\begin{equation}\label{2}
    \varphi(u+v)+\varphi(u-v)=2[\varphi(u)+\varphi(v)],
    \quad u,
    v\in Y.
\end{equation}
Denote by $\Gamma(X)$ the set of Gaussian distributions on
    $X$. We note that according this definition the generated distributions are Gaussian.

Denote by ${\rm Aut}(X)$ the group of topological automorphisms of
the group $X$, and by $I$ the identity
automorphism of a group. If  $\alpha\in{\rm Aut}(X)$,
then the adjoint  automorphism $\tilde\alpha\in{\rm Aut}(Y)$ is defined
as follows
   $(x, \tilde\alpha y) = (\alpha x, y)$ for all
 $x \in X$, $y \in Y$. Note that  $\alpha\in {\rm Aut}(X)$
 if and only if $\tilde\alpha\in {\rm Aut}(Y)$.
 Denote by $\mathbb{R}$   the group of real numbers,  by
$\mathbb{T}=\{z\in \mathbb{C}: |z|=1\}$ the circle group    (the one dimensional torus),
and by $\mathbb{Z}$ the group of integers.
Let $n$ be an integer.
Denote by $f_n$ the mapping of $X$ into $X$ defined by the formula $f_nx=nx$.
Put $X_{(n)}={\rm Ker} f_n$ and $X^{(n)}=f_n(X)$. A group $X$ is called
a Corwin group if $X^{(2)}=X$.

\bigskip

\noindent{\textbf{2.  Cram\'er's  theorem  for $Q$-independent random variables}}

\bigskip

\noindent In the article \cite{Fe1}, see also \cite[Theorem 4.6]{Fe5a}, the following analogue
of the classical Cram\'er theorem for the random variables with values in
a locally compact Abelian group was proved.

\medskip

\noindent{\bf Theorem А.}  {\it  Let $X$ be a locally compact Abelian group. Assume that $X$ contains no
 subgroup topologically isomorphic to the circle group
 $\mathbb{T}$. Let $Y$ be the character group of $X$.
 Let $\xi_1$ and $\xi_2$ be independent random variables with values
 in the group $X$. If a random variable
$\xi=\xi_1+\xi_2$ has a Gaussian distribution, then  $\xi_j$,
$j=1, 2,$ are also Gaussian. To put it in another way, if
 $\mu_{\xi}\in \Gamma(X)$ and $$\hat\mu_{\xi}(y)=\hat\mu_{\xi_1}(y)\hat\mu_{\xi_2}(y),$$
then $\hat\mu_{\xi_j}(y)$, $j=1, 2,$ are the characteristic functions
of Gaussian distributions}.

We prove that Theorem A remains true if
we change the condition of independence for $Q$-independence.
The following statement is valid.

\medskip

\noindent{\bf Proposition 1.} {\it Let $X$ be a locally compact Abelian group. Assume that $X$ contains no
 subgroup topologically isomorphic to the circle group
 $\mathbb{T}$. Let $\xi_1$ and $\xi_2$ be $Q$-independent random variables with values
 in the group $X$. If a random variable
$\xi=\xi_1+\xi_2$ has a Gaussian distribution, then  $\xi_j$,
$j=1, 2,$ are also Gaussian.}

\medskip

  To prove Proposition 1 we need the following lemma which is
a group analogue of the classical Marcinkiewicz theorem.

\medskip

\noindent{\bf Lemma 1} (\cite{Fe6}, see also \cite[Theorem 5.11]{Fe5a}).
{\it  Let $X$ be a locally compact Abelian group.  Assume that $X$ contains no
 subgroup topologically isomorphic to the circle group
 $\mathbb{T}$.  Let $Y$ be the character group of $X$ and
  $f(y)$ be a characteristic function on the group $Y$. If $f(y)$
is of the form
$$
f(y)=\exp\{P(y)\}, \quad  y\in Y,
$$
where $P(y)$  is a continuous polynomial,  then $P(y)$ is a polynomial of degree
 $\le 2$, and $f(y)$ is the characteristic function of a Gaussian distribution}.

\medskip

\noindent{\it Proof of Proposition} 1. Since the random variables
$\xi_1$ and $\xi_2$ are $Q$-independent,  we have
 \begin{equation}\label{3}
    \hat\mu_\xi(y)={\bf E}[(\xi, y)]={\bf E}[(\xi_1+\xi_2, y)]=
    {\bf E}[(\xi_1, y)(\xi_2, y)]=\hat\mu_{\xi_1}(y)\hat\mu_{\xi_2}(y)\exp\{q(y,  y)\},
    \quad y\in Y,
\end{equation}
 where $q(y_1,  y_2)$ is a continuous polynomial on the group
  $Y^2$. It follows from the definition of a polynomial on a group   that $q(y,  y)$ is a continuous polynomial on the group $Y$.
 By the condition
 the characteristic function $\hat\mu_\xi(y)$  is represented in the form (\ref{1}).
Then it follows from  (\ref{3}) that
\begin{equation}\label{4}
    (-x, y)\hat\mu_{\xi_1}(y)\hat\mu_{\xi_2}(y)=\exp\{P(y)\},
    \quad y\in Y,
\end{equation}
where $P(y)=-\varphi(y)-q(y,  y)$. Since $\varphi(y)$ is a continuous polynomial on the group
$Y$, $P(y)$ is also a continuous polynomial on the group
$Y$. The left-hand side of (\ref{4}) is a characteristic function.
Since the group $X$ contains no
 subgroup topologically isomorphic to the circle group $\mathbb{T}$, by Lemma 1
the left-hand side of (\ref{4}) is the characteristic function
of a Gaussian distribution. This implies by Theorem A  that the random variables
$\xi_j$, $j=1, 2,$ are also Gaussian. $\Box$

\medskip

We note that since independent random variables are  $Q$-independent,
and Theorem A fails if a locally compact Abelian group $X$ contains a subgroup
topologically isomorphic to the circle group $\mathbb{T}$ (see e.g.
\cite[\S 4]{Fe5a}), the condition on the group $X$ in Proposition 1 can not be weaken.

\bigskip

\noindent{\textbf{3. Skitovich--Darmois's  theorem
for $Q$-independent random variables}}

\bigskip

\noindent In the article \cite{Fe16}, see also  \cite[\S 10]{Fe5a},
the analogue
of the Skitovich--Darmois  theorem  for random variables with values in
a locally compact Abelian group was proved.

\medskip

\noindent{\bf Theorem B.} {\it Let $X$ be a locally compact Abelian group.  Assume that $X$ contains no
 subgroup topologically isomorphic to the circle group
 $\mathbb{T}$.  Let $\xi_1, \dots, \xi_n$ be   independent random variables with values
 in the group $X$ such that their characteristic functions do not vanish.
 Consider the linear forms
$L_1=\alpha_1\xi_1+\cdots+\alpha_n\xi_n$ and
$L_2=\beta_1\xi_1+\cdots+\beta_n\xi_n$, where coefficients
$\alpha_j,\beta_j\in {\rm Aut}(X)$. If the linear forms  $L_1$ and $L_2$
are independent, then the random variables $\xi_j$ are Gaussian.}

\medskip

We prove that Theorem B remains true if
we change the condition of independence of $\xi_1, \dots, \xi_n$ and
$L_1$, $L_2$ for $Q$-independence.
The following statement holds true.

\medskip

\noindent{\bf Theorem 1.} {\it Let $X$ be a locally compact Abelian group. Assume that $X$ contains no
 subgroup topologically isomorphic to the circle group
 $\mathbb{T}$.  Let $\xi_1, \dots, \xi_n$ be   $Q$-independent random variables with values
 in the group $X$ such that their characteristic functions do not vanish.
 Consider the linear forms
$L_1=\alpha_1\xi_1+\cdots+\alpha_n\xi_n$ and
$L_2=\beta_1\xi_1+\cdots+\beta_n\xi_n$, where coefficients
$\alpha_j,\beta_j\in {\rm Aut}(X)$. If the linear forms  $L_1$ and $L_2$
are $Q$-independent, then the random variables $\xi_j$ are Gaussian.}

\medskip

To prove Theorem 1 we need the following lemmas.

\medskip

\noindent{\bf Lemma 2}. {\it Let $X$ be a locally compact Abelian group and
$Y$ be its character group. Let  $\xi_1, \dots, \xi_n$ be
$Q$-independent random variables with values in the group
$X$, and let  $\alpha_j,\beta_j\in {\rm Aut}(X)$. The linear forms
$L_1=\alpha_1\xi_1+\cdots+\alpha_n\xi_n$ and
$L_2=\beta_1\xi_1+\cdots+\beta_n\xi_n$ are $Q$-independent if and only if the characteristic functions  $\hat\mu_{\xi_j}(y)$ satisfy the equation
\begin{equation}\label{13}
\prod_{j=1}^n\hat\mu_{\xi_j}(\tilde\alpha_ju+
\tilde\beta_jv)=\left(\prod_{j=1}^n\hat\mu_{\xi_j}(\tilde\alpha_ju)
\prod_{j=1}^n\hat\mu_{\xi_j}(\tilde\beta_jv)\right)\exp\{q(u, v)\},
\quad u, v\in Y,
\end{equation}
where $q(u, v)$ is a continuous polynomial on the group $Y^2$, $q(0, 0)=0$.}

\medskip

 \noindent{\it Proof}.
On the one hand, since the random variables $\xi_1,
\dots, \xi_n$  are $Q$-independent,
the join characteristic function of $L_1$ and $L_2$ is of the form
\begin{equation}\label{11}
\hat\mu_{(L_1, L_2)}(u, v)={\bf E}[(L_1, u)(L_2, v)]={\bf E}[(\alpha_1\xi_1+\cdots+\alpha_n\xi_n, u)(\beta_1\xi_1+\cdots+
\beta_n\xi_n, v)]=$$$$={\bf E}[(\xi_1, \tilde\alpha_1u+\tilde\beta_1v)\cdots(\xi_n, \tilde\alpha_nu+\tilde\beta_nv)]=$$$$=
\left(\prod_{j=1}^n\hat\mu_{\xi_j}(\tilde\alpha_ju+
\tilde\beta_jv)\right)\exp\{q_1(\tilde\alpha_1u+\tilde\beta_1v,\dots, \tilde\alpha_nu+\tilde\beta_nv)\},
\quad u, v\in Y,
\end{equation}
where $q_1(y_1, \dots, y_n)$ is a continuous polynomial on the group  $Y^n$.
On the other hand, since the random variables $\xi_1,
\dots, \xi_n$  are $Q$-independent, we have
 \begin{equation}\label{22_02_1}
\hat\mu_{L_1}(y)={\bf E}[(L_1, y)]={\bf E}[(\alpha_1\xi_1+\cdots+\alpha_n\xi_n, y)]=
\left(\prod_{j=1}^n\hat\mu_{\xi_j}(\tilde\alpha_jy)\right)
\exp\{q_1(\tilde\alpha_1y,\dots, \tilde\alpha_ny)\},
\quad y\in Y,
\end{equation}
\begin{equation}\label{22_02_2}
\hat\mu_{L_2}(y)={\bf E}[(L_2, y)]={\bf E}[(\beta_1\xi_1+\cdots+\beta_n\xi_n, y)]=
\left(\prod_{j=1}^n\hat\mu_{\xi_j}(\tilde\beta_jy)\right)
\exp\{q_1(\tilde\beta_1y,\dots, \tilde\beta_ny)\},
\quad y\in Y.
\end{equation}
Assume that the linear forms   $L_1$ and $L_2$ are $Q$-independent.
Then the join characteristic function of $L_1$ and $L_2$ can be written in
 the form
\begin{equation}\label{22_02_3}
\hat\mu_{(L_1, L_2)}(u, v)=\hat\mu_{L_1}(u)\hat\mu_{L_2}(v)
\exp\{q_2(u, v)\},
\quad u, v\in Y,
\end{equation}
where $q_2(u, v)$ is a continuous polynomial on the group $Y^2$.
Put
\begin{equation}\label{22_02_4}
q(u, v)=-q_1(\tilde\alpha_1u+\tilde\beta_1v,\dots, \tilde\alpha_nu+\tilde\beta_nv)+q_1(\tilde\alpha_1u,\dots, \tilde\alpha_nu)+q_1(\tilde\beta_1v,\dots, \tilde\beta_nv)+q_2(u, v).
\end{equation}
It follows from the definition of a polynomial on a group   that $q(u, v)$ is a continuous polynomial on the group $Y^2$. Obviously, (\ref{13}) follows from (\ref{11})--(\ref{22_02_4}).
If (\ref{13}) holds, then   (\ref{22_02_3}) follows from (\ref{11})--(\ref{22_02_2}),
where $q_2(u, v)$ is defined by formula (\ref{22_02_4}) and, obviously, $q_2(u, v)$ is a continuous polynomial on the group $Y^2$. So, lemma is proved.  It should be noted that the proof of the lemma, except
the fact that $q(u, v)$ is a polynomial, is the same as in the case when  $X$ is the group of real numbers (see  \cite{KS}). $\Box$

\medskip

\noindent{\bf Lemma 3}. {\it  Let $Y$ be an Abelian group and $b_j$ be automorphisms of
the group $Y$  such that $b_i\ne b_j$ for $i\ne j$. Consider
on the group
$Y$ the equation
\begin{equation}\label{5}
\sum_{j = 1}^{n}  \psi_j(u + b_j v ) = P(u) + Q(v) +R(u, v),
\quad u, v \in Y,
\end{equation}
where $\psi_j(y)$, $P(y)$ и $Q(y)$ are functions on
  $Y$, and $R(u, v)$ is a polynomial on $Y^2$. Then
$P(y)$ is a polynomial on $Y$.}

\medskip

\noindent{\it Proof}. In proving lemma we use the finite-difference method.
Let $h_1$ be an arbitrary element of the group
$Y$. Put $k_1=-{b^{-1}_n}h_1$. Then
$h_1+b_n k_1=0$. Substitute $u+h_1$ for $u$ and
$v+k_1$ for $v$  in equation (\ref{5}). Subtracting
equation (\ref{5}) from the resulting equation we obtain
\begin{equation}
\label{6}
    \sum_{j = 1}^{n-1} \Delta_{l_{1j}}{\psi_j(u + b_j v)}
    =\Delta_{h_{1}} P(u)+\Delta_{k_{1}} Q(v)+\Delta_{(h_{1},
    k_{1})}R(u, v),
\quad u,v\in Y,
\end{equation}
where $l_{1j}= h_1+b_j k_1=(b_j-b_n)k_1$, $j=1, 2,\dots, n-1$.  Let
$h_2$
be an arbitrary element of the group $Y$. Put
$k_2=-{b^{-1}_{n-1}} h_2$. Then $h_2+b_{n-1} k_2 =0$.
Substitute $u+h_2$ for $u$ and
$v+k_2$ for $v$  in equation (\ref{6}). Subtracting
equation (\ref{6}) from the resulting equation we find
\begin{equation}\label{7}
    \sum_{j = 1}^{n-2} \Delta_{l_{2j}} \Delta_{l_{1j}}{\psi_j(u + b_j v)}
    =\Delta_{h_{2}} \Delta_{h_{1}} P(u)+\Delta_{k_{2}} \Delta_{k_{1}} Q(v)+\Delta_{(h_{2}, k_{2})} \Delta_{(h_{1}, k_{1})}R(u, v),
\quad u, v\in Y,
\end{equation}
where $l_{2j}= h_2+b_j k_2=(b_j-b_{n-1})k_2$, $j=1, 2,\dots,n-2$.
Arguing   similarly  as above  we get the equation
\begin{equation}\label{8}
   \Delta_{l_{n-1,1}} \Delta_{l_{n-2,1}}\dots \Delta_{l_{11}}
   {\psi_1(u + b_1 v)}
    =\Delta_{h_{n-1}}\Delta_{h_{n-2}}\dots
    \Delta_{h_{1}} P(u)+$$$$+\Delta_{k_{n-1}} \Delta_{k_{n-2}}\dots \Delta_{k_{1}} Q(v)+\Delta_{(h_{n-1}, k_{n-1})} \Delta_{(h_{n-2}, k_{n-2})}\dots \Delta_{(h_{1}, k_{1})}R(u, v),
\quad u, v\in Y,
\end{equation}
where $h_m$ are   arbitrary elements of the group $Y$,
$k_m=-{b^{-1}_{n-m+1}}h_m$,
   $m=1, 2,\dots, n-1$, $l_{mj}=
h_m+b_j k_m=(b_j-b_{n-m+1})k_m$, $j=1, 2,\dots, n-m$. Let  $h_n$ be an arbitrary
element  of the group
 $Y$. Put $k_n=-{b^{-1}_1}h_n$.
Then $h_n+b_1 k_n=0$.
Substitute $u+h_n$ for $u$ and
$v+k_n$ for $v$  in equation (\ref{8}). Subtracting
equation (\ref{8}) from the resulting equation we obtain
\begin{equation}\label{9}
   \Delta_{h_{n}}\Delta_{h_{n-1}}\dots
    \Delta_{h_{1}} P(u)+\Delta_{k_{n}} \Delta_{k_{n-1}}
    \dots \Delta_{k_{1}} Q(v)+$$$$+\Delta_{(h_{n}, k_{n})}
    \Delta_{(h_{n-1}, k_{n-1})}\dots
    \Delta_{(h_{1}, k_{1})}R(u, v)=0,
\quad u, v\in Y.
\end{equation}
Let $h_{n+1}$ be an arbitrary
element  of the group $Y$.
Substitute $h_{n+1}$ for $u$   in equation (\ref{9}). Subtracting
equation (\ref{9}) from the resulting equation we find
\begin{equation}\label{10}
   \Delta_{h_{n+1}} \Delta_{h_{n}}\Delta_{h_{n-1}}\dots
    \Delta_{h_{1}} P(u)+\Delta_{(h_{n+1}, 0)}\Delta_{(h_{n}, k_{n})} \Delta_{(h_{n-1}, k_{n-1})}\dots \Delta_{(h_{1}, k_{1})}R(u, v)=0,
\quad u, v\in Y.
\end{equation}

We note that if  $h$ and $k$ are arbitrary elements
of the group $Y$, by the condition
\begin{equation}\label{z16}
  \Delta^{l+1}_{(h, k)}R(u, v)=0,
\quad u, v \in Y,
\end{equation}
for some   $l$. Since    $h_m$,
$m=1, 2, \dots,n+1$ are  arbitrary elements
of the group  $Y$, we can put in (\ref{10}) $h_1=\dots=h_{n+1}=h$,  and apply
to the both sides of the resulting equation the operator
$\Delta^{l+1}_{(h, k)}$. Taking into account (\ref{z16}), we get
$$
   \Delta_{h}^{l+n+2}P(u)=0, \quad u, h\in Y.
$$
So, the lemma is proved. $\Box$

\medskip

\noindent{\it Remark} 1. Let $l$ be the degree of the polynomial $R(u, v)$. Using some properties of
polynomials on groups (see e.g.   \cite[\S 5]{Fe5a}), it is not difficult
to obtain from (\ref{10}) that the degree of the polynomial
 $P(u)$ does not exceed $\max\{n, l\}$.

\medskip

\noindent{\it Proof of Theorem} 1. By Lemma  2 the characteristic functions
 $\hat\mu_{\xi_j}(y)$ satisfy equation (\ref{13}).
 Put  $g_j(y)=|\hat\mu_j(\tilde\alpha_j y)|^2$, $b_j=\tilde\alpha_j^{-1}\tilde\beta_j$, $R(u, v)=2{\rm Re}\ q(u, v)$.
 It follows from (\ref{13})  that
\begin{equation}\label{14}
 \prod_{j=1}^n g_j(u+
b_jv)=\left(\prod_{j=1}^ng_j(u)
\prod_{j=1}^ng_j(b_jv)\right)\exp\{R(u, v)\},
\quad u, v\in Y.
\end{equation}
Put $\psi_j(y)=\log g_j(y)$. Assume first that $b_i\ne b_j$ for $i\ne j$.
 It follows from (\ref{14}) that
\begin{equation}\label{15}
\sum_{j = 1}^{n}  \psi_j(u + b_j v ) = P(u) + Q(v) +
R(u, v),\quad u, v \in Y,
\end{equation}
where $$P(y)=\sum_{j=1}^n\psi_j(y), \quad Q(y)=
\sum_{j=1}^n\psi_j(b_jy).$$ By Lemma 3
$P(y)$ is a polynomial on $Y$. Obviously, the polynomial $P(y)$
is continuous. We have
\begin{equation}\label{16}
\prod_{j=1}^n g_j(y)=\exp\{P(y)\}, \quad y\in Y.
\end{equation}
Since the group $X$ contains no subgroup topologically isomorphic to the
circle group $\mathbb{T}$ and the left-hand side in (\ref{16})
is a characteristic function, by Lemma 1  the left-hand side in (\ref{16})
is the characteristic function of a Gaussian distribution.
This implies by Theorem A that all $g_j(y)$ are characteristic functions of   Gaussian distributions.
Applying Theorem A again we get that $\hat\mu_{\xi_j}(\tilde\alpha_jy)$ are characteristic functions of   Gaussian distributions, and hence, $\hat\mu_{\xi_j}(y)$ are characteristic functions of   Gaussian distributions, i.e.
the random variables $\xi_j$ are Gaussian.

Assume now that not all $b_j$ are different automorphisms.
Changing if  necessary functions  numbering   $g_j(y)$, we can assume that    $b_1=\dots =b_{l_1}$, \dots,
$b_{l_{k-1}+1}=\dots =b_{l_k}=b_n$. Put
$h_1(y)=g_1(y)\cdots g_{l_1}(y)$,
\dots, $h_k(y)=g_{l_{k-1}+1}(y)\cdots g_{l_k}(y)$. Reasoning as above we get that all  $h_l(y)$ are characteristic function of   Gaussian distributions.
By Theorem A this implies that all $g_j(y)$, and hence all  $\hat\mu_{\xi_j}(y)$ are characteristic functions of   Gaussian distributions, i.e.
the random variables $\xi_j$ are Gaussian.
$\Box$

\medskip

The following statement results from the proof of Theorem 1.

\medskip

\noindent{\bf Corollary 1}. {\it Let $X$ be a locally compact Abelian group. Assume that $X$ contains no
 subgroup topologically isomorphic to the circle group
 $\mathbb{T}$. Let $Y$ be the character group of $X$.
 Let $\xi_1, \dots, \xi_n$ be independent random variables with values
 in the group $X$ such that their characteristic functions do not vanish
 and satisfy equation
 $(\ref{13})$, where
$\tilde\alpha_j,\tilde\beta_j\in {\rm Aut}(Y)$. Then
$\hat\mu_{\xi_j}(y)$ are characteristic functions of   Gaussian distributions.}

\medskip

\noindent{\it Remark} 2. It should be noted that since independent random variables are $Q$-independent and in the proof of Theorem 1 we do not use
Theorem B, Theorem B follows from Theorem 1.

\medskip

\noindent{\it Remark} 3. It is well known, see e.g.  \cite[Lemma 7.8]{Fe5a},
if a locally compact Abelian group $X$ contains
a subgroup topologically isomorphic to the circle group
 $\mathbb{T}$, then there exist independent random variables $\xi_1$ and $\xi_2$ with values
 in the group $X$ such that their characteristic functions do not vanish,
$\xi_1+\xi_2$ and $\xi_1-\xi_2$ are independent, and
$\mu_{\xi_1}, \mu_{\xi_2}\notin \Gamma(X)$. Since independent random variables are $Q$-independent, the condition on the group
$X$ in Theorem 1 can not be weaken.

\bigskip

\noindent{\textbf{4.  Heyde's  theorem
for $Q$-independent random variables}}

\bigskip

\noindent The following  group analogue of the
 well-known Heyde theorem on characterization of
 a Gaussian distribution on the real line (\cite{He}, see also
\cite[\S 13.4]{KaLiRa}), was proved in \cite{Fe30}.

\medskip

\noindent{\bf Theorem С.}  {\it  Let  $X$ be a locally compact Abelian group.  Assume that $X$ contains no elements of order $2$.  Let $\alpha$ be a topological automorphism of the group $X$.
Let  $\xi_1$ and  $\xi_2$ be independent random variables with values in
$X$ with non-vanishing characteristic functions. The symmetry of the conditional distribution of the linear form
$L_2 = \xi_1 + \alpha\xi_2$ given  $L_1 = \xi_1 +
\xi_2$ implies that   $\xi_1$ and $\xi_2$ are Gaussian  if and only if
  $\alpha$ satisfies the  condition
   \begin{equation}\label{z1}
{\rm Ker}(I+\alpha)=\{0\}.
\end{equation}}

\medskip

We prove that Theorem C remains true if
we change the condition of independence of random variables $\xi_1$ and $\xi_2$ for $Q$-independence.
The following statement is valid.

\medskip

\noindent{\bf Theorem 2.}  {\it  Let  $X$ be a locally compact Abelian group.  Assume that $X$ contains no elements of order $2$.  Let $\alpha$ be a topological automorphism of the group $X$.
Let  $\xi_1$ and  $\xi_2$ be $Q$-independent random variables with values in
$X$ with non-vanishing characteristic functions. The symmetry of the conditional distribution of the linear form
$L_2 = \xi_1 + \alpha\xi_2$ given  $L_1 = \xi_1 +
\xi_2$ implies that   $\xi_1$ and $\xi_2$ are Gaussian  if and only if
  $\alpha$ satisfies  condition $(\ref{z1})$}.

\medskip

To prove Theorem 2 we need the following lemma.

\medskip

\noindent {\bf  Lemma 4}. {\it  Let  $X$ be a locally compact Abelian group, $Y$ be its character group, and $\alpha$ be a topological automorphism of the group $X$.
Let  $\xi_1$ and  $\xi_2$ be $Q$-independent random variables with values in
$X$. If the conditional distribution of the linear form
$L_2 = \xi_1 + \alpha\xi_2$ given  $L_1 = \xi_1 +
\xi_2$ is symmetric, then the characteristic functions
 $\hat\mu_{\xi_j}(y)$  satisfy the equation
\begin{equation}\label{z2}
\hat\mu_{\xi_1}(u+v )\hat\mu_{\xi_2}(u+\tilde\alpha v )=
\hat\mu_{\xi_1}(u-v )\hat\mu_{\xi_2}(u-\tilde\alpha v)\exp\{q(u, v)\}, \quad u, v \in Y,
\end{equation}
where $q(u, v)$ is a continuous polynomial on the group $Y^2$, $q(0, 0)=0$.}

\medskip

\noindent{\it Proof}. Since the random variables $\xi_1$ and  $\xi_2$ are $Q$-independent,   the join characteristic function of $L_1$ and $L_2$ is of the form
\begin{equation}\label{z3}
\hat\mu_{(L_1, L_2)}(u, v)={\bf E}[(L_1, u)(L_2, v)]={\bf E}[(\xi_1 + \xi_2, u)(\xi_1 + \alpha\xi_2, v)]=$$$$={\bf E}[(\xi_1, u+ v) (\xi_2, u+\tilde\alpha v)]=
\hat\mu_{\xi_1}(u+v )\hat\mu_{\xi_2}(u+\tilde\alpha v)\exp\{q_1(u+v, u+\tilde\alpha v)\},
\quad u, v\in Y,
\end{equation}
where $q_1(u, v)$ is a continuous polynomial on the group $Y^2$, $q_1(0, 0)=0$.
Similarly, the join characteristic function of $L_1$ and $-L_2$ is of the form
\begin{equation}\label{z4}
\hat\mu_{(L_1, -L_2)}(u, v)=
\hat\mu_{\xi_1}(u-v )\hat\mu_{\xi_2}(u-\tilde\alpha v)\exp\{q_1(u-v, u-\tilde\alpha v)\},
\quad u, v\in Y.
\end{equation}
The symmetry of the conditional distribution of the linear form
$L_2$ given  $L_1$ means that the random vectors
$(L_1, L_2)$ and $(L_1, -L_2)$ are identically distributed, i.e.
the join characteristic functions (\ref{z3}) and (\ref{z4})
are the same. Put $q(u, v)=-q_1(u+v, u+\tilde\alpha v)+q_1(u-v, u-\tilde\alpha v)$. It is obvious that $q(u, v)$ is a continuous polynomial on the group $Y^2$, and (\ref{z2}) follows from (\ref{z3})
and (\ref{z4}). $\Box$

\medskip

\noindent{\it Proof of Theorem} 2. Since   independent random variables are $Q$-independent, the necessity follows from Theorem C.

Sufficiency. We will use some ideas from the articles   \cite{Fe30}  and  \cite{My1}.
By Lemma 4 the characteristic functions
$\hat\mu_{\xi_j}(y)$  satisfy equation (\ref{z2}). Put $b=\tilde\alpha$.
Let $w, z\in Y$.  Substituting in (\ref{z2}) $u=b w$, $v=-w$ and $u=z$, $v=-z$, we get
\begin{equation}\label{z5}
\hat\mu_{\xi_1}((b - I)w)=\hat\mu_{\xi_1}((I+
b )w)\hat\mu_{\xi_2}(2b w)\exp\{q(b w, -w)\},
\quad w\in Y.
\end{equation}
\begin{equation}\label{z6}
\hat\mu_{\xi_2}(-(b - I)z)=\hat\mu_{\xi_1}(2z)\hat\mu_{\xi_2}((I+b) z)\exp\{q(z, -z)\},
\quad z\in Y.
\end{equation}
Substituting in  (\ref{z2}) $u=bw+z$, $v=w+z$, we obtain
\begin{equation}\label{z7}
\hat\mu_{\xi_1}((I+b)w +2z)\hat\mu_{\xi_2}(2b w +(I+b)z)$$$$=\hat\mu_{\xi_1}((b-I)w)
\hat\mu_{\xi_2}(-(b-I)z)\exp\{q(w+b z, w+z)\},
\quad w, z\in Y.
\end{equation}
It follows from (\ref{z5})--(\ref{z7}) that
\begin{equation}\label{z8}
\hat\mu_{\xi_1}((I+b)w +2z)\hat\mu_{\xi_2}(2b w +(I+b)z)=\hat\mu_{\xi_1}((I+
b )w)\hat\mu_{\xi_2}(2b w)\exp\{q(b w, -w)\}\times$$$$\hat\mu_{\xi_1}(2z)\hat\mu_{\xi_2}((I+b) z)\exp\{q(z, -z)\}\exp\{q(bw+z, w+z)\},
\quad w, z\in Y.
\end{equation}
Put $p(u, v)=q(b u, -u)+q(v, -v)+q(bu+v, u+v)$.
Obviously,  $p(u, v)$ is a continuous polynomial on the group $Y^2$, $p(0, 0)=0$.

Thus, the characteristic functions $\hat\mu_{\xi_j}(y)$  satisfy the equation
\begin{equation}\label{z9}
\hat\mu_{\xi_1}((I+b)u +2v)\hat\mu_{\xi_2}(2b u +(I+b)v)=$$$$=\hat\mu_{\xi_1}((I+
b )u)\hat\mu_{\xi_2}(2b u)\hat\mu_{\xi_1}(2v)\hat\mu_{\xi_2}((I+b) v)\exp\{p(u, v)\},
\quad u, v\in Y.
\end{equation}

Since the group $X$ contains no elements of order 2, $X$ contains
no subgroup topologically isomorphic to the circle group $\mathbb{T}$. In the case  when $I+\alpha$ and $f_2$ are topological automorphisms of the group $X$, the statement of the theorem follows from   (\ref{z9}) and Corollary 1. In general, $I+\alpha$ and $f_2$ are only continuous monomorphisms, and we reason as follows.

Put   $g_j(y)=|\hat\mu_j(y)|^2$, $j=1, 2$,  $R(u, v)=2{\rm Re}\ p(u, v)$.
 It follows from (\ref{z9}) that the functions $g_j(y)$  satisfy the equation
 \begin{equation}\label{z10}
g_1 ((I+b)u +2v)g_2 (2b u +(I+b)v)=$$$$=g_1 ((I+
b )u)g_2 (2b u)g_1 (2v) g_2 ((I+b) v)\exp\{R(u, v)\},
\quad u, v\in Y.
\end{equation}
Put $\psi_j(y)= \log g_j(y)$, $j=1, 2$.  Then
$(\ref{z10})$ implies that the functions $\psi_j(y)$ satisfy the equation
\begin{equation}\label{z11}
\psi_1((I+b) u+2 v)+\psi_2(2b  u+(I+b) v)=P(u)+Q(v)+R(u, v), \quad u, v \in
Y,
\end{equation}
where
\begin{equation}\label{z12}
P(y)=\psi_1((I+b)
y)+\psi_2(2b
y), \quad Q(y)=\psi_1(2 y)+\psi_2((I+b) y).
\end{equation}
 In solving the equation (\ref{z11}) we use the finite-difference method.
Let $h_1$ be an arbitrary element of the group
$Y$. Substitute in
 (\ref{z11}) $u+(I+b) h_1$ for $u$ and $v-2 b  h_1$ for $v$. Subtracting
equation
  (\ref{z11}) from the resulting equation we obtain
  \begin{equation}\label{z13}
    \Delta_{(I-b)^2 h_1}{\psi_1((I+b) u+2 v)}
    =\Delta_{(I+b) h_1} P(u)+\Delta_{-2b  h_1} Q(v)+\Delta_{((I+b) h_1, -2b  h_1)}R(u, v),
\quad u, v\in Y.
\end{equation}
Let $h_2$ be an arbitrary element of the group
$Y$. Substitute in
(\ref{z13})   $u+2h_{2}$ for $u$ and
$v-(I+b)h_{2}$ for $v$. Subtracting
equation
  (\ref{z13}) from the resulting equation we get
 \begin{equation}\label{z14}
     \Delta_{2 h_2}\Delta_{(I+b) h_1} P(u)+\Delta_{-(I+b) h_2}\Delta_{-2b  h_1} Q(v)+$$$$+\Delta_{(2h_{2}, -(I+b)h_{2})}\Delta_{((I+b) h_1, -2b  h_1)}R(u, v)=0,
\quad u, v\in Y.
\end{equation}
Let $h_3$ be an arbitrary element of the group
$Y$. Substitute in
(\ref{z14})   $u+h_3$ for $u$. Subtracting
equation
  (\ref{z14}) from the resulting equation we find
 \begin{equation}\label{z15}
   \Delta_{h_3}\Delta_{2 h_2}\Delta_{(I+b) h_1} P(u)+\Delta_{(h_3, 0)}\Delta_{(2h_{2}, -(I+b)h_{2})}\Delta_{((I+b) h_1, -2b  h_1)}R(u, v)=0,
\quad u, v\in Y.
\end{equation}

We note that if $h$ and $k$ are arbitrary elements of the group
 $Y$, then (\ref{z16}) holds. Applying
to the both sides of  equation (\ref{z15}) the operator
$\Delta^{l+1}_{(h, k)}$ we obtain
 \begin{equation}\label{24_02_1}
   \Delta^{l+1}_{(h, k)}\Delta_{h_3}\Delta_{2 h_2}\Delta_{(I+b) h_1} P(u)=0,
\quad u, v\in Y.
\end{equation}

 Since the group $X$ contains no elements of order 2,
the subgroup $Y^{(2)}$ is dense in  $Y$ (\cite[$($24.22)]{HeRo1}). Since
$\alpha$ satisfies condition
$(\ref{z1})$, the subgroup
$(I+b)(Y)$ is also dense in
$Y$ (\cite[(24.41)]{HeRo1}). Taking into account that
$h_1, h_2, h_3$ are arbitrary elements of the group $Y$,   (\ref{z16}) and
(\ref{24_02_1}) imply  that
$$
   \Delta_{h}^{l+4}P(u)=0, \quad u, h\in Y,
$$
i.e. $P(y)$ is a continuous polynomial. It follows from (\ref{z12})
that $\exp\{P(y)\}$ is a characteristic function.
Since the group $X$ contains no elements of order 2,    $X$ contains
no subgroup topologically isomorphic to the circle group $\mathbb{T}$.
 Then, by Lemma 1 $P(y)$ is a continuous polynomial of degree $\le 2$, and
  $\exp\{P(y)\}=g_1((I+b)y)g_2(2by)$ is the characteristic function of a Gaussian distribution. By Theorem A, this implies that $g_1((I+b)y)$ and $g_2(2by)$ are
  the characteristic functions of Gaussian distributions.
 Since the subgroups
$(I+b)(Y)$ and $Y^{(2)}$   are dense in $Y$,  $g_1(y)$ and $g_2(y)$
  are also
  the characteristic functions of Gaussian distributions.
  Applying Theorem A again we get that $\xi_1$ and $\xi_2$ are Gaussian. $\Box$

\medskip

\noindent{\it Remark} 4. Let  $X$ be a locally compact Abelian group, $Y$ be its character group, and $\alpha$ be a topological automorphism of the group $X$.
Let  $\xi_1$ and  $\xi_2$ be  random variables with values in
$X$. Consider the linear forms
$L_2 = \xi_1 + \alpha\xi_2$ and  $L_1 = \xi_1 +
\xi_2$. Note that the   symmetry of the conditional distribution of the linear form
$L_2$ given  $L_1$  means that the random vectors
$(L_1, L_2)$ and $(L_1, -L_2)$ are identically distributed. Let $\eta_1$ and  $\eta_2$ be random variables
with values in the group
 $X$.
 Following A. Kagan and G. Sz\'ekely  (\cite{KS}) we say that   $\eta_1$ and  $\eta_2$
 are
$Q$-identically distributed if
\begin{equation}\label{xy1}
\hat\mu_{\eta_1}(y)=\hat\mu_{\eta_2}(y)\exp\{q(y)\}, \quad y\in Y,
\end{equation}
where $q(y)$ is a continuous polynomial on the group $Y$.
We will also assume that $q(0)=0$.
 It is easy to see that Lemma 4
remains true if we change in the lemma the condition of
 symmetry of the conditional distribution of the linear form
$L_2$ given  $L_1$ for the condition that the random vectors
$(L_1, L_2)$ and $(L_1, -L_2)$, i.e. the random variables  with values in the group $X^2$,
 are $Q$-identically distributed.
Since  in the proof of Theorem 2 we only use equation (\ref{z2}), Theorem 2 remains true if
  we change in Theorem 2 the condition of
 symmetry of the conditional distribution of the linear form
$L_2$ given  $L_1$ for weaker condition: the random vectors
$(L_1, L_2)$ and $(L_1, -L_2)$ are $Q$-identically distributed.

\bigskip

\noindent{\textbf{5.  Кас-Bernstein's  theorem
for $Q$-independent random variables}}

\bigskip

\noindent To prove the main theorem of this section we need some facts on structure of locally compact Abelian groups
and duality theory (see \cite[Chapter 6]{HeRo1}).   Any locally compact Abelian group
$X$ is topologically isomorphic to a group of the form $ \mathbb{R}^m\times G,$
where $m\ge 0$, and the group $G$ contains a compact open subgroup.
Let $Y$ be a character group of the group $X$.
We will also denote the character group of the group $X$ by $X^*$.
If $K$ is a closed subgroup of $X$, denote by
 $A(Y, K) = \{y \in Y: (x, y) = 1$ \mbox{ for all } $x \in K \}$
its annihilator.
Note that $A(X, A(Y, K))=K$, and the following topological automorphisms
 $K^*\cong Y/A(Y, K)$ and
$(X/K)^*\cong A(Y, K)$ hold.
Moreover, if $B$ is a closed subgroup of $Y$, then
  any character of the group $B$ is of the form $y\mapsto (x, y)$ for some
  $x\in X$. We also note that $A\left(Y, X_{(n)}\right)=
\overline{Y^{(n)}}$ for any natural $n$.

  Let  $\mu\in{\rm M}^1(X)$. Denote by  $\sigma(\mu)$
 the support of $\mu$.
Denote by $E_x$  the degenerate distribution
 concentrated at an element $x\in X$.
If $K$ is a compact subgroup of the group $X$, then denote by $m_K$ the Haar distribution of
the group $K$.
We note that the characteristic function of a distribution
$m_K$ is of the form
\begin{equation}\label{30}
\hat m_K(y)=
\begin{cases}
1, & \text{\ if\ }\   y\in A(Y, K),
\\  0, & \text{\ if\ }\ y\not\in
A(Y, K).
\end{cases}
\end{equation}
Denote by $I_B(X)$  the set of Haar distributions $m_K$
of compact Corwin subgroups  $K$ of the group $X.$

In the article \cite{Fe5}, see also  \cite[Theorem 7.10]{Fe5a},
the following analogue of the Кас-Bernstein  theorem
for  random variables with values in a locally compact Abelian group was proved.

\medskip

\noindent{\bf Theorem D.} {\it Let $X$ be a locally compact Abelian group.
Assume that the connected component of zero of the group $X$
 contains no elements of order $2$. Let   $\xi_1$ and $\xi_2$ be
 independent random variables
with values in  $X$  such
that $\xi_1+\xi_2$ and $\xi_1-\xi_2$ are independent.
Then $\mu_{\xi_1},
\mu_{\xi_2}\in \Gamma(X)*I_B(X)$. Moreover,
\begin{equation}\label{51}
\mu_{\xi_1}=\mu_{\xi_2}*E_x,
\end{equation}
where $x\in X$.}

\medskip

We note that in contrast to Theorems B and C, in Theorem D we do not assume
that the characteristic functions of considering random variables do not vanish.
We will prove that Theorem D, except, generally speaking, statement   (\ref{51}),
 remains true if we substitute the condition of independence for
  $Q$-independence. The following statement holds true.

\medskip

\noindent{\bf Theorem 3.} {\it Let $X$ be a locally compact Abelian group.
Assume that the connected component of zero of the group $X$
 contains no elements of order $2$. Let   $\xi_1$ and $\xi_2$ be
 $Q$-independent random variables
with values in  $X$  such
that $\xi_1+\xi_2$ and $\xi_1-\xi_2$ are $Q$-independent.
Then   $\mu_{\xi_j}=\gamma_j*m_W$, where
$\gamma\in \Gamma(X)$, $j=1, 2$, and $W$ is a compact Corwin subgroup.}

\medskip

To prove Theorem 3 we need the following lemmas.

\medskip

\noindent{\bf Lemma 5} (see e.g. \cite[Proposition 5.7]{Fe5a}).
{\it
Let $Y$ be a compact Abelian group and $f(y)$ be a continuous polynomial on
$Y$. Then $f(y)\equiv const$}.

\medskip

We formulate as a lemma the following well-known statement.

\medskip

\noindent{\bf Lemma 6}. {\it Let $X$ be a locally compact Abelian group and
 $Y$ be its character group. Let  $\mu\in{\rm
M}^1(X)$. Put $E=\{y\in Y:\ \hat\mu(y)=1\}$. Then
$\sigma(\mu)\subset A(X,E)$.}

\medskip

\noindent{\bf Lemma 7} (\cite[\S 7]{Fe5a}). {\it
Let $X$ be a locally compact Abelian group of the form $X=\mathbb{R}^m\times
K,$ where $m\ge 0$, and $K$ is a compact Corwin group. Assume that the connected component of zero of the group $X$
 contains no elements of order $2$.
Let $Y$ be the character group of  $X$.
Then $X$ and $Y$ are groups with unique division by $2$.}

\medskip

\noindent{\it Proof of Theorem} 3. We will use the scheme of the proof of Theorem 7.10 in \cite{Fe5a}.  Taking into account the structure theorem for locally compact Abelian groups, we can assume without loss of generality, that  $X=\mathbb{R}^m\times
G,$ $Y=\mathbb{R}^m\times H$, where $m \ge 0$, $H\cong G^*$, and each of
 the groups  $G$ and $H$ contains a compact open subgroup.
Denote by $L$  a compact open subgroup in $H$. Put
$$
N_1=\{y\in Y: \hat\mu_{\xi_1}(y)\ne 0\}, \quad N_2=\{y\in Y:
\hat\mu_{\xi_2}(y)\ne 0\}, \quad N=N_1\cap N_2.
$$
By Lemma 2, the characteristic functions
$\hat\mu_{\xi_j}(y)$ satisfy eqaution (\ref{13}) which takes
the form
\begin{equation}\label{18}
\hat\mu_{\xi_1}(u+v)\hat\mu_{\xi_2}(u-v)=
 \hat\mu_{\xi_1}(u)\hat\mu_{\xi_2}(u)\hat\mu_{\xi_1}(v)
 \hat\mu_{\xi_2}(-v)\exp\{q(u, v)\}, \quad u, v\in Y,
\end{equation}
where $q(u, v)$ is a continuous polynomial on the group $Y^2$, $q(0, 0)=0$. It follows from
(\ref{18}) that $N$ is a subgroup of $Y$. Obviously,   $N$
is an open subgroup. Consider the intersection  $B=N\cap L$. Since
any open subgroup is closed, $B$ is a compact open subgroup of $H$.
By Lemma 5 $q(u, v)=0$, $u, v\in B$. This implies that the restriction
of equation (\ref{18}) to the subgroup $B$ is of the form
 \begin{equation}\label{19}
\hat\mu_{\xi_1}(u+v)\hat\mu_{\xi_2}(u-v)=
 \hat\mu_{\xi_1}(u)\hat\mu_{\xi_2}(u)\hat\mu_{\xi_1}(v)
 \hat\mu_{\xi_2}(-v), \quad u, v\in B.
\end{equation}
Put $A=B^*$. Since $B$ is a compact group,  $A$ is a discrete group, and hence, $A$ is totally disconnected. In particular,
the connected component of zero of the group $A$
 contains no elements of order $2$. Taking into account that
 Gaussian distributions on an arbitrary locally compact totally disconnected Abelian group are degenerated  (see  \cite[Chapter IV]{Pa}), and
   the characteristic functions $\hat\mu_{\xi_j}(y)$ do not vanish on $B$,
    by Theorem D applying the group A we obtain
$$
\hat\mu_{\xi_1}(y)=(x_1, y), \quad \hat\mu_{\xi_2}(y)=(x_2, y),
\quad y\in B,
$$
where $x_j\in X$, $j=1, 2$. Consider the new random variables   $\xi_j'=\xi_j-x_j$.
It is obvious that the random variables $\xi_j'$ are $Q$-independent,
and $\xi_1'+\xi_2'$ and
$\xi_1'-\xi_2'$ are also $Q$-independent. Therefore, passing from the random variables
$\xi_j$ to the random variables $\xi_j'$, we can prove the theorem assuming from the beginning that
$$
\hat\mu_{\xi_1}(y)=\hat\mu_{\xi_2}(y)=1, \quad y\in B,
$$
is fulfilled. Then, by Lemma 6 $\sigma(\mu_{\xi_j})\subset A(X, B)$, $j=1, 2$.
We have $A(X,
B)\cong(Y/B)^*=((\mathbb{R}^m\times
H)/B)^*\cong\mathbb{R}^m\times(H/B)^*$. Put $F=(H/B)^*$.
Since $B$ is an open subgroup of $H$, the factor-group  $H/B$
is discrete, and hence, $F$ is a compact group.

Thus, we reduced the proof of the theorem to the case, when
$X=\mathbb{R}^m\times F$, $Y=\mathbb{R}^m\times D$, where $F$ is a compact group, and   $D\cong F^*$ is a discrete group.
Let $D_2$  be the subgroup of $D$ consisting of all elements
$y\in D$ such that the order of $y$ is a power of 2.  It follows from Lemma 5 that $q(u, v)=0$,
$u, v\in D_2$. Hence, the restriction of equation (\ref{18})
to the subgroup $D_2$ is of the form
 (\ref{19}).
 Substituting in equation (\ref{19}) $u=v=y$    and  $u=-v=y$, we obtain
\begin{equation}
\label{20}
\hat\mu_{\xi_1}(2y)=(\hat\mu_{\xi_1}(y))^2|\hat\mu_{\xi_2}(y)|^2,
\quad\hat\mu_{\xi_2}(2y)
=|\hat\mu_{\xi_1}(y)|^2(\hat\mu_{\xi_2}(y))^2, \quad y \in D_2.
\end{equation}
It follows from (\ref{20}) that for any natural $n$ the following identities
\begin{equation}
\label{21}
|\hat\mu_{\xi_j}(2^ny)|=|\hat\mu_{\xi_1}(y)\hat\mu_{\xi_2}(y)|^{2^{2n-1}},
 \quad y \in D_2, \quad j=1, 2,
\end{equation}
hold  true. Let $y\in D_2$. Then $2^ny=0$ for some natural $n$, and (\ref{21}) implies that
$|\hat\mu_{\xi_1}(y)|=|\hat\mu_{\xi_2}(y)|=1$ for $y\in D_2$.
It follows from this that there exist elements    $x'_j\in X$, $j=1, 2,$
such that
$$
\hat\mu_{\xi_1}(y)=(x'_1, y), \quad \hat\mu_{\xi_2}(y)=(x'_2, y),
\quad y\in D_2.
$$
Reasoning as above we can prove the theorem assuming that the  identities
$$
\hat\mu_{\xi_1}(y)=\hat\mu_{\xi_2}(y)=1, \quad y\in D_2,
$$
are fulfilled from the beginning. Then, by Lemma 6  $\sigma(\mu_{\xi_j})\subset A(X, D_2)$, $j=1, 2$.
Put $M=A(X, D_2)$. It is obvious that $M= \mathbb{R}^m\times K$,
where $K$ is a compact group. We have $M^*\cong Y/D_2$. It is clear that
the factor-group
$Y/D_2$ contains no elements of order $2$.
This implies that
 $\overline {M^{(2)}}=M$. Since
$\overline{M^{(2)}}=M^{(2)}$, it means that $M$ is a Corwin  group, and hence, $K$ is also a Corwin  group.

Thus, we reduced the proof of the theorem to the case, when
$X=\mathbb{R}^m\times K$,  where $K$ is   a compact Corwin  group.
Since the connected component of zero of the group $X$
 contains no elements of order $2$, by Lemma 7
$X$ and $Y$ are   groups  with unique division by $2$. We will check that this implies the equality $N_1=N_2$.

Substituting in equation (\ref{18}) $u=v=y$   and  $u=-v=y$, we get
\begin{equation}
\label{22}
\hat\mu_{\xi_1}(2y)=(\hat\mu_{\xi_1}(y))^2|
\hat\mu_{\xi_2}(y)|^2\exp\{q(y,y)\},\quad\hat\mu_{\xi_2}(2y)
=|\hat\mu_{\xi_1}(y)|^2(\hat\mu_{\xi_2}(y))^2\exp\{q(y,-y)\}, \quad y \in Y.
\end{equation}
Assume for definiteness that there exists an element $y_0\in Y$ such that $\hat\mu_{\xi_1}(y_0)\ne 0$,
$\hat\mu_{\xi_2}(y_0)=0$. Since $Y$ is a group  with unique division by $2$, we have $y_0=2y'$.
Substituting $y=y'$ in (\ref{22}), we get from the one hand,  $\hat\mu_{\xi_1}(y')\hat\mu_{\xi_2}(y')\ne 0$, and from the other hand, either $\hat\mu_{\xi_1}(y')= 0$ or $\hat\mu_{\xi_2}(y')= 0$. This contradiction implies that $N_1=N_2=N$.
Since $Y$ is a group  with unique division by $2$,   (\ref{22}) implies that   $N$ is also  a group  with unique division by $2$.
Put $W=A(X, N)$. It is not difficult to check that $W$ is a compact Corwin group. Since  $(X/W)^*\cong N$, this implies that
 $X/W$  is   a group  with unique division by $2$. Obviously, the group
 $X/W$ contains no subgroup topologically isomorphic to the circle group  $\mathbb{T}$.
 Consider the restriction of equation (\ref{18}) to the group $N$ and apply
 Corollary 1 to the group
$X/W$. We obtain the representation
\begin{equation}
\label{23}
\hat\mu_{\xi_1}(y)=(x_1, y)\exp\{-\varphi_1(y)\},\quad\hat\mu_{\xi_2}(y)=(x_2, y)\exp\{-\varphi_2(y)\},  \quad y \in N,
\end{equation}
where $x_j\in X$, and  $\varphi_j(y)$ are continuous non-negative  functions on  $N$  satisfying equation (\ref{2}).
Extend the functions $\varphi_j(y)$ from the subgroup
$N$ to the group $Y$ in such a way that the extended functions are continuous  non-negative and also satisfy equation (\ref{2}) (see e.g. \cite[Lemma 3.18]{Fe5a}). We retain the notation  $\varphi_j(y)$ for the extended functions.  Let
$\gamma_j$ be the Gaussian distributions on the group
 $X$ with the characteristic functions
\begin{equation}
\label{24} \hat\gamma_j(y)=(x_j,
y)\exp\{-\varphi_j(y)\}, \quad y\in Y.
\end{equation}
Since $N=A(Y, W)$, it follows from (\ref{30})
that the characteristic function of the Haar distribution
 $m_W$ is of the form
\begin{equation}
\label{25} \hat m_W(y)=
\begin{cases}
1, & \text{\ if\ }\ \   y\in N,
\\ 0,& \text{\ if\ }\ \ y\not\in N.
\end{cases}
\end{equation}
It follows from (\ref{23})--(\ref{25}) that
$\hat\mu_{\xi_j}(y)=\hat\gamma_j(y)\hat m_W(y)$.
Hence,  $\mu_{\xi_j}=\gamma_j*m_W$, $j=1, 2$. $\Box$

\medskip

\noindent{\it Remark} 5. Let $X$ be a locally compact Abelian group and
 $Y$ be its character group.
Assume that the connected component of zero of $X$ is non-zero and contains no elements
of order  $2$.
Then, in contrast to Theorem  D, we can not state in Theorem 3  that
 $\mu_{\xi_1}=\mu_{\xi_2}*E_x$ for some  $x\in X.$  Indeed, since the connected component of zero of $X$ is non-zero, there exists a non-degenerate Gaussian distribution on $X$  (see  \cite[Chapter IV]{Pa}). Let
$\gamma_j$ be a non-degenerate Gaussian distributions on the group $X$  with the characteristic functions
\begin{equation}\label{31}
\hat\gamma_j(y)=\exp\{-\varphi_j(y)\}, \quad y\in Y, \quad j=1, 2,
\end{equation}
where $\varphi_j(y)$ are continuous non-negative functions on $Y$, satisfying equation (\ref{2}).
Let $\xi_1$ and $\xi_2$ be independent random variables with values in the group $X$ such that $\mu_{\xi_j}=\gamma_j$.
Put $\eta_1=\xi_1+\xi_2$, $\eta_2=\xi_1-\xi_2$. On the one hand,
taking into account that $\varphi_j(-y)=\varphi_j(y)$ and independence of  the random variables $\xi_j$, the characteristic functions of the random variables
$\eta_j$ are of the form
\begin{equation}\label{w1}
\hat\mu_{\eta_j}(y)={\bf E}[(\eta_1, y)]=\hat\gamma_1(y)\hat\gamma_2(y),  \quad y\in Y, \quad j=1, 2.
\end{equation}
On the other hand, taking into account independence of  the random variables  $\xi_j$ and (\ref{w1}), the join characteristic function of the random variables $\eta_1$ and $\eta_2$
can be written in the form
$$
\hat\mu_{(\eta_1, \eta_2)}(u, v)={\bf E}[(\eta_1, u)(\eta_2, v)]={\bf E}[(\xi_1+\xi_2, u)(\xi_1-\xi_2, v)]={\bf E}[(\xi_1, u+v)(\xi_2,u -v) ]=$$$$=\hat\gamma_1(u+v)\hat\gamma_2(u-v)=\hat\mu_{\eta_1}(u)
\hat\mu_{\eta_2}(v)
\exp\{-\varphi_1(u+v)-\varphi_2(u-v)+\varphi_1(u)+\varphi_2(u)+
\varphi_1(v)+\varphi_2(v)\}.
$$
Put
$$
q(u, v)=-\varphi_1(u+v)-\varphi_2(u-v)+\varphi_1(u)+\varphi_2(u)+
\varphi_1(v)+\varphi_2(v).
$$
Then $q(u, v)$ is a continuous polynomial on $Y^2$. It is easy to see that $q(u, v)\equiv 0$
if and only if
$\varphi_1(y)\equiv\varphi_2(y)$. Thus,  if $\varphi_1(y)\not\equiv\varphi_2(y)$, then
the random variables $\eta_1$ and $\eta_2$ are not independent, but $Q$-independent.

From what has been said it follows that if  $\varphi_1(y)\not\equiv\varphi_2(y)$, then
the random variables   $\xi_1$ and $\xi_2$ are $Q$-independent, they even are
independent,  the random variables $\xi_1+\xi_2$ and $\xi_1-\xi_2$ are $Q$-independent, i.e. the conditions of Theorem 3 are fulfilled, whereas the distributions  $\mu_{\xi_1}=\gamma_1$ and $\mu_{\xi_2}=\gamma_2$ are not shifts one another.

\medskip

\noindent{\it Remark} 6. It is interesting to observe that there exist
$Q$-independent random variables with values in a locally compact Abelian
group $X$ such that they are not independent if and only if
 the connected component of zero of the group $X$ is non-zero.
  Indeed, let $X$ be a totally disconnected group and $Y$ be its character group.
  Then all elements of the group $Y$  are compact. Assume that  $\xi_1, \xi_2, \dots, \xi_n$ are
$Q$-independent random variables with values in $X$. By Lemma 5 any continuous polynomial
 on $Y$ is a constant on the subgroup of all compact elements of the group
   $Y$. Hence,   in (\ref{17}) $q(y_1, \dots, y_n)\equiv 0.$  It means that the random variables
     $\xi_1, \xi_2, \dots, \xi_n$
 are independent.

Assume that the connected component of zero of a group $X$ is non-zero.
Let $\xi_1$ and $\xi_2$ be independent non-degenerated Gaussian random variables with values in the group
$X$ such that $\mu_{\xi_1}=\gamma_1$ and $\mu_{\xi_2}=\gamma_2$, where the characteristic
functions  of the distributions $\gamma_j$ are of the form (\ref{31}). If
 $\varphi_1(y)\not\equiv\varphi_2(y)$, then as noted in Remark 5,
the random variables   $\eta_1=\xi_1+\xi_2$ and $\eta_2=\xi_1-\xi_2$ are
$Q$-independent, but not independent.

\medskip

\noindent{\it Remark} 7. Compare Theorems  1, 2 and 3 with Theorems  B, C and D.
We see that Theorems  B, C and D remain true for the corresponding groups if
we change the condition of independence for $Q$-independence.
 It turns out that either a characterization theorem remains true or not after such change
 depends
 on the group. We give an example of a characterization theorem
 on the circle group $\mathbb{T}$ which fails if  we change the condition of independence for $Q$-independence. Moreover,  actually a weak analogue of this theorem fails too.

The following theorem results from the article \cite{BaESta}, see also  \cite[Theorem 9.9]{Fe5a}. It characterizes a Gaussian distribution on the circle group  ${\mathbb T}$ (compare with
 Theorem B).

\medskip

\noindent{\bf Theorem E.}
{\it  Let $\xi_1$ and $\xi_2$ be independent identically distributed
random variables with values in the circle group  ${\mathbb T}$ such that their characteristic functions do not vanish.   Assume that $\xi_1+\xi_2$ and
$\xi_1-\xi_2$ are independent. Then  $\mu_{\xi_j}\in \Gamma({\mathbb T})$,
$j=1, 2$.}

\medskip

The character group of the circle group ${\mathbb T}$ is topologically isomorphic to
  ${\mathbb Z}$. We will assume that ${\mathbb T}^*={\mathbb Z}$. We construct
  independent identically distributed
random variables  $\xi_1$ and $\xi_2$ with values in the circle group ${\mathbb T}$  and with distribution
   $\mu_{\xi_j}=\gamma$,  $j=1, 2$, such that  the random variables $\xi_1+\xi_2$ and $\xi_1-\xi_2$ are   $Q$-independent, and the characteristic function
$\hat\gamma(n)$ is represented in the form
$$
\hat\gamma(n)=\exp\{-\varphi(n)\}, \quad n\in {\mathbb Z},
$$
where $\varphi(n)$
is a polynomial on the group of integers ${\mathbb Z}$ of sufficiently arbitrary form.

  Let $\varphi(n)$ be a polynomial on the group of integers ${\mathbb Z}$ such that  $\varphi(0)=0$, $\varphi(-n)=\varphi(n)$, $n\in{\mathbb Z}$, and
$$
\sum_{n\in \mathbb{Z}} \exp\{-\varphi(n)\}< 2.
$$
This implies that
$$
\rho(t)=\sum_{n\in \mathbb{Z}} \exp\{-\varphi(n)-int\}>0,
\quad t\in\mathbb{R}.
$$
It is also obvious that
$$
{1\over 2\pi}\int\limits_{-\pi}^\pi\rho(t)dt=1.
$$
Let $\gamma$  be a distribution on the circle group  ${\mathbb T}$ with the density $
r(e^{it})=\rho(t)$ with respect to $m_{\mathbb T}$. Then
$\hat\gamma(n)=\exp\{-\varphi(n)\}$, $n\in {\mathbb Z}$.

Let $\xi_1$ and  $\xi_2$ be independent identically distributed
random variables with values in the circle group ${\mathbb T}$  and with distribution
$\mu_{\xi_j}=\gamma$.  Reasoning as in Remark 5, it is easy to make sure that $\xi_1+\xi_2$ and $\xi_1-\xi_2$ are $Q$-independent.

\end{document}